\documentclass[a4paper, BCOR=0.0mm, DIV=calc]{scrartcl}
\usepackage[T1]{fontenc}
\usepackage[utf8]{inputenc}
\usepackage[ngerman]{babel}
\usepackage[osf,sc]{mathpazo}
\usepackage{textcomp}
\usepackage{microtype}
\usepackage{amsmath, amssymb, amsfonts}
\usepackage{tabularx}
\usepackage{nicefrac}
\KOMAoptions{twoside=false, twocolumn=false, headinclude=false, footinclude=false, mpinclude=false, pagesize=auto}
\recalctypearea
\newcommand{\ebinom}[2]{\left(\frac{#1}{#2} \right)}
\setkomafont{section}{\mdseries\scshape\Large}
\setkomafont{title}{\mdseries\scshape}

\begin{document}
\title{Beobachtung über Reihen, deren Terme nach den Sinus oder Kosinus vielfacher Winkel fortschreiten\footnote{
Originaltitel: "`Observationes generales circa series, quarum termini secundum sinus vel cosinus angulorum multiplorum progrediuntur"', erstmals publiziert in "`\textit{Nova Acta Academiae Scientarum Imperialis Petropolitinae} 7, 1793, pp. 87-98 "', Nachdruck in "`\textit{Opera Omnia}: Series 1, Volume 16, pp. 163 - 177"', Eneström-Nummer E655, übersetzt von: Alexander Aycock, Textsatz: Artur Diener,  im Rahmen des Projektes "`Eulerkreis Mainz"' }}
\author{Leonhard Euler}
\date{}
\maketitle
\paragraph{§1}
Wenn also die Summation dieser Reihe
\[
	A + Bx + Cxx + Dx^3 + \mathrm{etc}
\]
bekannt war, sodass, welcher Wert auch immer dem Buchstaben $x$ zugeteilt wird, ihre Summe angegeben werden kann, dann wird auch immer so die Summe dieser Reihe
\[
	A + B\cos{\varphi} + C\cos{2\varphi} + D\cos{3\varphi} + \mathrm{etc}
\]
wie dieser
\[
	B\sin{\varphi} + C\sin{2\varphi} + D\sin{3\varphi} + E\sin{4\varphi} + \mathrm{etc}
\]
beschafft werden können. Weil nämlich die Summe der ersten Reihe durch eine gewisse Funktion von $x$ ausgedrückt wird, die wir mit dem Charakter  $\Delta : x$ bezeichnen wollen, sodass
\[
	\Delta : x = A + Bx + Cxx + Dx^3 + \mathrm{etc}
\]
ist, wird, wenn wir anstelle von $x$
\[
	\cos{\varphi} + \sqrt{-1}\sin{\varphi}
\]
sowie
\[
	\cos{\varphi} - \sqrt{-1}\sin{\varphi}
\]
schreiben, die Summen der daher entstehenden Reihen
\[
	2A + 2B\cos{\varphi} + 2C\cos{2\varphi} + 2D\cos{3\varphi} + 2E\cos{4\varphi} + \mathrm{etc}
\]
sein, deren Summe also
\[
	\Delta : \left( \cos{\varphi} + \sqrt{-1}\sin{\varphi} \right) + \Delta : \left( \cos{\varphi} - \sqrt{-1}\sin{\varphi} \right)
\]
sein wird; wenn wir aber letztere von der ersten abziehen, wird diese Reihe hervorgehen:
\[
	2B\sqrt{-1}\sin{\varphi} + 2C\sqrt{-1}\sin{2\varphi} + 2D\sqrt{-1}\sin{3\varphi} + 2E\sqrt{-1}\sin{4\varphi} + \mathrm{etc}
\]
hervorgehen, deren Summe also
\[
	\Delta : \left( \cos{\varphi} + \sqrt{-1}\sin{\varphi} \right) - \Delta : \left( \cos{\varphi} - \sqrt{-1}\sin{\varphi} \right)
\]
sein wird.
\paragraph{§2}
Damit wir diese Ausdrücke vereinfachen, wollen wir der Kürze wegen
\[
	\cos{\varphi} + \sqrt{-1}\sin{\varphi} = p \quad \text{und} \quad \cos{\varphi} - \sqrt{-1}\sin{\varphi} = q
\]
setzen und es wird, wie im Allgemeinen bekannt ist,
\[
	pq = 1
\]
sein und daher
\[
	q = \frac{1}{p};
\]
dann wird aber
\[
	\cos{\varphi} = \frac{p+q}{2},\quad \cos{2\varphi} = \frac{pp+qq}{2},\quad \cos{3\varphi} = \frac{p^3 + q^3}{2},\quad \cos{4\varphi} = \frac{p^4 + q^4}{2},\quad \mathrm{etc}
\]
sein. Außerdem aber wird man für die Sinus haben
\[
	\sin{\varphi} = \frac{p-q}{2\sqrt{-1}},\quad \sin{2\varphi} = \frac{pp-qq}{2\sqrt{-1}},\quad \sin{3\varphi} = \frac{p^3 - q^3}{2\sqrt{-1}},\quad \sin{4\varphi} = \frac{p^4 - q^4}{2\sqrt{-1}},\quad \mathrm{etc},
\]
nach Festsetzen wovon wir diese $2$ Summationen erreichen:
\[
	A\cos{0\varphi} + B\cos{\varphi} + C\cos{2\varphi} + D\cos{3\varphi} + \mathrm{etc} = \frac{\Delta : p + \Delta : q}{2}
\]
und
\[
	A\sin{0\varphi} + B\sin{\varphi} + C\sin{2\varphi} + D\sin{3\varphi} + \mathrm{etc} = \frac{\Delta : p - \Delta : q}{2\sqrt{-1}}.
\]
\paragraph{§3}
Wir wollen nun für die anfängliche Reihe eine beliebige entwickelte Potenz des Binoms nehmen, welche
\[
	(1+x)^n = 1 + \frac{n}{1}x + \frac{n(n-1)}{1\cdot 2}xx + \frac{n(n-1)(n-2)}{1\cdot 2\cdot 3}x^3 + \mathrm{etc}
\]
ist, sodass in diesem Fall
\[
	\Delta : x = (1+x)^n
\]
ist, dann aber, damit wir diesen Ausdruck zusammenfassen, wollen wir die einzelnen Koeffizienten, wie ich es schon des öfteren gemacht habe, mit diesen Charakteren bezeichnen:
\[
	\ebinom{n}{0}, \ebinom{n}{1}, \ebinom{n}{2}, \ebinom{n}{3}, \ebinom{n}{4}, \mathrm{etc},
\]
sodass
\begin{align*}
	\ebinom{n}{0} &= 1 \\
	\ebinom{n}{1} &= n \\
	\ebinom{n}{2} &= \frac{n(n-1)}{1\cdot 2} \\
	\ebinom{n}{3} &= \frac{n(n-1)(n-2)}{1\cdot 2\cdot 3} \\
	\mathrm{etc}
\end{align*}
ist, wo es förderlich sein wird, bemerkt zu haben, dass im Allgemeinen
\[
	\ebinom{n}{1} = \ebinom{n}{n-i}
\]
ist und daher
\[
	\ebinom{n}{n} = \ebinom{n}{0} = 1.
\]
Außerdem ist in der Tat klar, dass, sooft $i$ entweder eine negative Zahl war oder eine positive größer als $n$, dass dann immer
\[
	\ebinom{n}{i} = 0
\]
ist, wenn natürlich $n$ eine ganze Zahl war. Nachdem diese Dinge also bemerkt worden sind, werden wir diese anfängliche Summation haben:
\[
	(1+x)^n = \ebinom{n}{0} + \ebinom{n}{1}x + \ebinom{n}{2}x^2 + \ebinom{n}{3}x^3 + \ebinom{n}{4}x^4 + \mathrm{etc},
\]
woher also durch die gerade angegebenen Festlegungen diese beiden Summationen berechnen werden:
\[
	\ebinom{n}{0}\cos{0\varphi} + \ebinom{n}{1}\cos{1\varphi} + \ebinom{n}{2}\cos{2\varphi} + \ebinom{n}{3}\cos{3\varphi} + \mathrm{etc} = \frac{(1+p)^n + (1+q)^n}{2}
\]
und
\[
	\ebinom{n}{0}\sin{0\varphi} + \ebinom{n}{1}\sin{1\varphi} + \ebinom{n}{2}\sin{2\varphi} + \ebinom{n}{3}\sin{3\varphi} + \mathrm{etc} = \frac{(1+p)^n - (1+q)^n}{2\sqrt{-1}}.
\]
In diesem Fall aber werden sich, obwohl die für $p$ und $q$ angenommenen Formeln imaginär sind, dennoch die Formeln auf reelle Werte zurückführen lassen, so wie wir in den folgenden Problemen zeigen werden.
\section*{Problem 1}
Nachdem diese Reihe der Kosinus vorgelegt wurde
\[
	1 + \frac{n}{1}\cos{\varphi} + \frac{n}{1}\cdot\frac{n-1}{2}\cos{2\varphi} + \frac{n}{1}\cdot\frac{n-1}{2}\cdot\frac{n-2}{3}\cos{3\varphi} + \mathrm{etc} = s,
\]
sodass durch die festgesetzten Charaktere
\[
	s = \ebinom{n}{0}\cos{0\varphi} + \ebinom{n}{1}\cos{1\varphi} + \ebinom{n}{2}\cos{2\varphi} + \ebinom{n}{3}\cos{3\varphi} + \mathrm{etc},
\]
ist ihre Summe reell auszudrücken.
\section*{Lösung}
\paragraph{§4}
Weil also, wie wir gerade gesehen haben
\[
	s = \frac{(1+p)^n + (1+q)^n}{2}
\]
ist, während
\[
	p = \cos{\varphi} + \sqrt{-1}\sin{\varphi}\quad \text{und} \quad q = \cos{\varphi} - \sqrt{-1}\sin{\varphi}
\]
wird, geht die ganze Aufgabe darauf zurück, dass dieser für $s$ beschaffte Ausdruck vom Imaginären befreit wird; es ist nämlich klar, wenn die Formeln $(1+p)^n$ und $(1+q)^n$ entwickelt werden, dass sich dann die imaginären Anteile von selbst aufheben werden, weil ja die zu summierende Reihe selbst entsteht; deshalb werden wir eine andere Auflösung suchen, dass ohne angewandte Entwicklung die imaginäre Anteile aufgehoben werden; das wird auf die folgende Weise gemacht werden können.
\paragraph{§5}
Weil $pq=1$ ist, wird die Formel $1+p$ so ausgedrückt werden können, dass
\[
	1+p = \left( \sqrt{p} + \sqrt{q} \right)\sqrt{p}
\]
ist; und auf ähnliche Weise wird
\[
	1+q = \left( \sqrt{p} + \sqrt{q} \right) \sqrt{q}
\]
sein; nach der Einführung dieser Werte wird unsere Summe als
\[
	s = \frac{1}{2} \left( \sqrt{p} + \sqrt{q} \right)^n \left( p^{\frac{n}{2}} + q^{\frac{n}{2}} \right)
\]
hervorgehen. Weil schon im Allgemeinen
\[
	p^{\alpha} + q^{\alpha} = 2\cos{\alpha\varphi}
\]
ist, wird
\[
	p^{\frac{1}{2}} + q^{\frac{1}{2}} = 2\cos{\tfrac{1}{2}\varphi} \quad \text{und} \quad p^{\frac{n}{2}} + q^{\frac{n}{2}} = 2\cos{\tfrac{1}{2}n\varphi}
\]
sein, nach Einsetzen welcher Werte die gesuchte  schon reell auf die folgende Weise ausgedrückt werden wird:
\[
	s = 2^n\cos^n{\tfrac{\varphi}{2}}\cos{\tfrac{1}{2}n\varphi}.
\]
\paragraph{§6}
Nach dieser Übereinkunft also haben wir eine besonders bemerkenswerte Summation erhalten, die sich so verhält, dass immer
\begin{align*}
	&1 + \frac{n}{1}\cos{\varphi} + \frac{n}{1}\cdot\frac{n-1}{2}\cos{2\varphi} + \frac{n}{1}\cdot\frac{n-1}{2}\cdot\frac{n-2}{3}\cos{3\varphi} + \mathrm{etc} \\
	&= 2^n\cos^n{\frac{1}{2}\varphi}\cos{\frac{1}{2}n\varphi},
\end{align*}
die immer mit der Wahrheit verträglich ist, welche Zahlen auch immer für $n$ eingesetzt werden, ob ganze oder gebrochene oder sogar negative. Es wird also der Mühe Wert sein, aus jeder Art einfachere Fälle vor Augen zu führen.
\section*{Entwicklung der Fälle, in denen der Exponent $n$ eine ganze positive Zahl ist}
\paragraph{§7}
Wir wollen die folgenden Fälle betrachten:
\begin{enumerate}
	\item Es sei $n=0$ und die Reihe selbst wird zur Einheit verschmelzen, die Summe wird aber gleich $1$ sein.
	\item Es sei $n=1$ und eine Reihe wird in 
	\[
		1 + \cos{\varphi}
	\]
	übergehen; die gefundene Summe aber liefert
	\[
		2\cos^2{\tfrac{1}{2}\varphi},
	\]
	es ist aber bekannt, dass
	\[
		2\cos^2{\tfrac{1}{2}\varphi} = 1 + \cos{\varphi}
	\]
	ist.
	\item Es sei $n=2$ und die Reihe wird in
	\[
		1 + 2\cos{\varphi} + \cos{2\varphi}
	\]
	übergehen, es entsteht aber die Summe
	\[
		4\cos^2{\tfrac{1}{2}\varphi}\cos{\varphi}.
	\]
	Gerade aber haben wir gesehen, dass $2\cos^2{\frac{1}{2}\varphi} = 1 + \cos{\varphi}$ ist, welche Form mit $2\cos{\varphi}$ multipliziert
	\[
		2\cos{\varphi}\cdot 2\cos^2{\tfrac{1}{2}\varphi} = 1 + 2\cos{\varphi} + \cos{2\varphi}
	\]
	ergibt.
	\item Es sei nun $n=3$ und es entsteht diese Reihe
	\[
		1 + 3\cos{\varphi} + 3\cos{2\varphi} + \cos{3\varphi},
	\]
	deren Summe gleich
	\[
		8\cos^3{\tfrac{1}{2}\varphi}\cos{\tfrac{3}{2}\varphi}
	\]
	ist, welche Formel durch hinreichend bekannte Reduktion die Reihe selbst ergibt.
	\item Es sei nun $n=4$ und die Reihe geht über in 
	\[
		1 + 4\cos{\varphi} + 6\cos{2\varphi} + 4\cos{3\varphi} + \cos{4\varphi},
	\]
	deren Summe also
	\[
		2^4\cos^4{\tfrac{1}{2}\varphi}\cos{2\varphi}
	\]
	sein wird, deren Gültigkeit man auch nicht schwer zeigt. Und so wird sich die Gültigkeit immer durch bekannte Reduktionen zeigen lassen. 
\end{enumerate}
\section*{Entwicklung der Fälle, in denen für $n$ eine ganze negative Zahl angenommen wird}
\paragraph{§8}
Wir wollen zuerst $n=-1$ setzen und es wird die folgende unendliche Reihe entstehen:
\[
	1 - \cos{\varphi} + \cos{2\varphi} - \cos{3\varphi} + \cos{4\varphi} - \cos{5\varphi} + \mathrm{etc}
\]
ins Unendlich, deren Summe also durch unsere allgemeine Reihe
\[
	\frac{\cos{\frac{1}{2}\varphi}}{2\cos{\frac{1}{2}\varphi}} = \frac{1}{2}
\]
sein wird, was freilich schon längst von den Mathematikern beobachtet worden ist. Wenn daher nämlich diese Reihe, deren Summe solange gleich $s$ gesetzt wird, mit $2\cos{\frac{1}{2}\varphi}$ multipliziert wird, wird man durch allbekannte Reduktionen
\[
	2s\cos{\frac{1}{2}\varphi} = \left. \begin{cases}
	\quad ~ 2\cos{\frac{1}{2}\varphi} - \cos{\frac{3}{2}\varphi} + \cos{\frac{5}{2}\varphi} - \cos{\frac{7}{2}\varphi} + \cos{\frac{9}{2}\varphi} \\
	\quad -\cos{\frac{1}{2}\varphi} + \cos{\frac{3}{2}\varphi} - \cos{\frac{5}{2}\varphi} + \cos{\frac{7}{2}\varphi} - \cos{\frac{9}{2}\varphi}
	\end{cases}\right\}
	\mathrm{etc}
\]
finden, was natürlich auf $2s\cos{\frac{1}{2}\varphi} = \cos{\frac{1}{2}\varphi}$ hinausläuft und daher $s = \frac{1}{2}$.
\paragraph{§9}
Wir wollen nun $n=-2$ setzen und es wird die folgende Reihe entstehen
\[
	1 - 2\cos{\varphi} + 3\cos{2\varphi} - 4\cos{3\varphi} + 5\cos{4\varphi} - 6\cos{5\varphi} + \mathrm{etc},
\]
deren Summe also gleich
\[
	\frac{\cos{\varphi}}{4\cos^2{\frac{1}{2}\varphi}}
\]
sein wird, deren Gültigkeit darüber hinaus auf die folgende Weise gezeigt werden kann. Nachdem die Summe der Reihe gleich $s$ gesetzt wurde, wird
\[
	2s\cos{\frac{1}{2}\varphi} = \left. \begin{cases}
	\quad ~~~ 2\cos{\frac{1}{2}\varphi} - 2\cos{\frac{3}{2}\varphi} + 3\cos{\frac{5}{2}\varphi} - 4\cos{\frac{7}{2}\varphi} \\
	\quad -2\cos{\frac{1}{2}\varphi} + 3\cos{\frac{3}{2}\varphi} - 4\cos{\frac{5}{2}\varphi} + 5\cos{\frac{7}{2}\varphi} 
	\end{cases} \right\}
	\mathrm{etc}
\]
sein, welcher Wert zur folgenden Reihe verschmilzt
\[
	2s\cos{\tfrac{1}{2}\varphi} = \cos{\tfrac{3}{2}\varphi} - \cos{\tfrac{5}{2}\varphi} + \cos{\tfrac{7}{2}\varphi} - \cos{\tfrac{9}{2}\varphi} + \mathrm{etc.}
\]
Man multipliziere erneut mit $2\cos{\frac{1}{2}\varphi}$ und es wird
\[
	4s\cos^2{\frac{1}{2}\varphi} = \left. \begin{cases}
	\quad \cos{\varphi} + \cos{2\varphi} - \cos{3\varphi} + \cos{4\varphi} - \cos{5\varphi} \\
	\quad \phantom{\cos{\varphi}} - \cos{2\varphi} + \cos{3\varphi} - \cos{4\varphi} + \cos{5\varphi}
	\end{cases}\right\}
	\mathrm{etc} = \cos{\varphi}
\]
hervorgehen und daher
\[
	s = \frac{\cos{\varphi}}{4\cos^2{\frac{1}{2}\varphi}},
\]
wie wir gefunden haben, oder es wird auch
\[
	s = \frac{\cos{\varphi}}{2(1+\cos{\varphi})}
\]
sein.
\paragraph{§10}
Es sei nun $n=-3$ und es wird diese unendliche Reihe entstehen
\[
	1 - 3\cos{\varphi} + 6\cos{2\varphi} - 10\cos{3\varphi} + 15\cos{4\varphi} - 21\cos{5\varphi} + \mathrm{etc},
\]
deren Summe also gleich
\[
	\frac{\cos{\frac{3}{2}\varphi}}{8\cos^3{\frac{1}{2}\varphi}}
\]
sein wird. Dieser Ausdruck wird aber weiter auf diesen zurückgeführt
\[
	s = \frac{1}{2} - \frac{3}{8\cos^2{\frac{1}{2}\varphi}} = \frac{1}{2} - \frac{3}{4(1+\cos{\varphi})},
\]
sodass auch
\[
	s = \frac{-1+2\cos{\varphi}}{4(1+\cos{\varphi})}
\]
ist.
\paragraph{§11}
Auf ähnliche Weise werden wir auch die folgenden Summationen erhalten:
\begin{align*}
	1 - 4\cos{\varphi} + 10\cos{2\varphi} - 20\cos{3\varphi} + 35\cos{4\varphi} - \mathrm{etc} &= \frac{\cos{2\varphi}}{16\cos^{4}{\frac{1}{2}\varphi}} \\
	1 - 5\cos{\varphi} + 15\cos{2\varphi} - 35\cos{3\varphi} + 70\cos{4\varphi} - \mathrm{etc} &= \frac{\cos{\frac{5}{2}\varphi}}{32\cos^{5}{\frac{1}{2}\varphi}} \\
	1 - 6\cos{\varphi} + 21\cos{2\varphi} - 56\cos{3\varphi} + 126\cos{4\varphi} - \mathrm{etc} &= \frac{\cos{3\varphi}}{64\cos^{6}{\frac{1}{2}\varphi}} \\
	1 - 7\cos{\varphi} + 28\cos{2\varphi} - 84\cos{3\varphi} + 210\cos{4\varphi} - \mathrm{etc} &= \frac{\cos{\frac{7}{2}\varphi}}{128\cos^{7}{\frac{1}{2}\varphi}} \\
	\mathrm{etc.}
\end{align*}
\section*{Entwicklung des Falles, in dem $n=\frac{1}{2}$ ist}
\paragraph{§12}
Es wird daher also die folgende unendliche Reihe gebildet werden
\[
	1 + \frac{1}{2}\cos{\varphi} - \frac{1\cdot 1}{2\cdot 4}\cos{2\varphi} + \frac{1\cdot 1\cdot 3}{2\cdot 4\cdot 6}\cos{3\varphi} - \frac{1\cdot 1\cdot 3\cdot 5}{2\cdot 4\cdot 6\cdot 8}\cos{4\varphi} + \mathrm{etc},
\]
deren Summe also gleich
\[
	\cos{\frac{1}{4}\varphi\sqrt{2\cos{\frac{1}{2}\varphi}}}
\]
sein wird, deren Gültigkeit nicht leicht sein wird anderswoher zu prüfen; in anderen Fällen aber springt sie natürlich ins Auge. Wenn z.\,B. $\varphi = 0$ war, wird man
\[
	1 + \frac{1}{2} - \frac{1\cdot 1}{2\cdot 4} + \frac{1\cdot 1\cdot 3}{2\cdot 4\cdot 6} - \frac{1\cdot 1\cdot 3\cdot 5}{2\cdot 4\cdot 6\cdot 8} + \mathrm{etc} = \sqrt{2}
\]
haben; die Reihe entsteht natürlich aus der Entwicklung 
\[
	(1+1)^{\frac{1}{2}} = \sqrt{2}.
\]
Wir wollen nun $\varphi = 180^{\circ}$ setzen, dass $\frac{1}{2}\varphi = 90^{\circ}$ ist, und die Reihe wird
\[
	1 - \frac{1}{2} - \frac{1\cdot 1}{2\cdot 4} - \frac{1\cdot 1\cdot 3}{2\cdot 4\cdot 6} - \frac{1\cdot 1\cdot 3\cdot 5}{2\cdot 4\cdot 6\cdot 8} - \mathrm{etc} = 0
\]
sein, was auch daher klar ist, weil diese Reihe aus der Form $(1-1)^{\frac{1}{2}}$ entsteht. Es sei auch $\varphi = 90^{\circ}$ und die daher entstehende Reihe wird
\[
	1 + \frac{1\cdot 1}{2\cdot 4} - \frac{1\cdot 1\cdot 3\cdot 5}{2\cdot 4\cdot 6\cdot 8} + \frac{1\cdot 1\cdot 3\cdot 5\cdot 7\cdot 9 }{2\cdot 4\cdot 6\cdot 8\cdot 10\cdot 12} - \mathrm{etc} = \cos{22^{\circ}\, 30'\sqrt[4]{2}}
\]
sein. Es ist aber
\[
	\cos{22^{\circ}\, 30'} = \sqrt{\frac{1+\cos{45^{\circ}}}{2}} = \sqrt{\frac{1+\sqrt{2}}{2\sqrt{2}}},
\]
woher man die Summe
\[
	\sqrt{\frac{1+\sqrt{2}}{\sqrt{2}}}
\]
folgert und so hat man diese höchst bemerkenswerte Summation
\[
	1 + \frac{1\cdot 1}{2\cdot 4} - \frac{1\cdot 1\cdot 3\cdot 5}{2\cdot 4\cdot 6\cdot 8} + \frac{1\cdot 1\cdot 3\cdot 5\cdot 7\cdot 9}{2\cdot 4\cdot 6\cdot 8\cdot 10\cdot 12} - \mathrm{etc} = \sqrt{\frac{1+\sqrt{2}}{2}}.
\]
Wir wollen auch $\varphi = 60^{\circ}$ nehmen und es wird diese Reihe entstehen
\begin{align*}
	&1 + \frac{1}{2}\cdot \frac{1}{2} + \frac{1\cdot 1}{2\cdot 4}\cdot\frac{1}{2} - \frac{1\cdot 1\cdot 3}{2\cdot 4\cdot 6}\cdot 1 + \frac{1\cdot 1\cdot 3\cdot 5}{2\cdot 4\cdot 6\cdot 8}\cdot\frac{1}{2} \\
	&+ \frac{1\cdot 1\cdot 3\cdot 5\cdot 7}{2\cdot 4\cdot 6\cdot 8\cdot 10}\cdot\frac{1}{2} - \frac{1\cdot 1\cdot 3\cdot 5\cdot 7\cdot 9}{2\cdot 4\cdot 6\cdot 8\cdot 10\cdot 12}\cdot 1 + \mathrm{etc},
\end{align*}
die Summe welcher Reihe also $\cos{15^{\circ}\sqrt[4]{3}}$ sein wird. Weil also
\[
	\cos{15^{\circ}} = \sqrt{\frac{1+\cos{30^{\circ}}}{2}} = \sqrt{\frac{2+\sqrt{3}}{4}}
\]
ist, wird die Summe der Reihe
\[
	\frac{1}{2}\sqrt{3+2\sqrt{3}}
\]
sein.
\section*{Entwicklung des Falles, in dem $n=-\frac{1}{2}$ ist}
\paragraph{§13}
Daher wird also die folgende unendliche Reihe gebildet werden
\[
	1 - \frac{1}{2}\cos{\varphi} + \frac{1\cdot 3}{2\cdot 4}\cos{2\varphi} - \frac{1\cdot 3\cdot 5}{2\cdot 4\cdot 6}\cos{3\varphi} + \frac{1\cdot 3\cdot 5\cdot 7}{2\cdot 4\cdot 6\cdot 8}\cos{4\varphi} - \mathrm{etc},
\]
deren Summe also gleich
\[
	\frac{\cos{\frac{1}{4}\varphi}}{\sqrt{2\cos{\frac{1}{2}\varphi}}}
\]
sein wird. Daher entsteht, wenn $\varphi = 0$ war, diese Summation
\[
	1 - \frac{1}{2} + \frac{1\cdot 3}{2\cdot 4} - \frac{1\cdot 3\cdot 5}{2\cdot 4\cdot 6} + \frac{1\cdot 3\cdot 5\cdot 7}{2\cdot 4\cdot 6\cdot 8} - \mathrm{etc} = \frac{1}{\sqrt{2}}.
\]
Diese Reihe entsteht nämlich aus der Form $(1+1)^{-\frac{1}{2}}$. Es sei nun $\varphi = 180^{\circ}$ und die resultierende Reihe wird
\[
	1 + \frac{1}{2} + \frac{1\cdot 3}{2\cdot 4} + \frac{1\cdot 3\cdot 5}{2\cdot 4\cdot 6} + \frac{1\cdot 3\cdot 5\cdot 7}{2\cdot 4\cdot 6\cdot 8} + \mathrm{etc} = \infty
\]
sein. Diese Reihe entsteht nämlich aus der Entwicklung $(1-1)^{-\frac{1}{2}}$. Wir wollen auch $\varphi = 90^{\circ}$ nehmen und die Reihe wird
\[
	1 - \frac{1\cdot 3}{2\cdot 4} + \frac{1\cdot 3\cdot 5\cdot 7}{2\cdot 4\cdot 6\cdot 8} - \frac{1\cdot 3\cdot 5\cdot 7\cdot 9\cdot 11}{2\cdot 4\cdot 6\cdot 8\cdot 10\cdot 12} + \mathrm{etc} = \frac{\cos{22^{\circ}\, 30'}}{\sqrt[4]{2}}
\]
sein. Vorher aber haben wir gesehen, dass
\[
	\cos{22^{\circ}\, 30'} = \sqrt{\frac{1+\sqrt{2}}{2\sqrt{2}}}
\]
ist, woher die Summe gleich
\[
	\frac{1}{2}\sqrt{1+\sqrt{2}}
\]
sein wird.
\paragraph{§14}
Im Allgemeinen wird es auch für beliebige Exponenten $n$ der Mühe Wert sein, dem Winkel $\varphi$ bestimmte Werte zuzuteilen; und nachdem freilich zuerst $\varphi = 0$ genommen wurde, werden wir
\[
	1 + \ebinom{n}{1} + \ebinom{n}{2} + \ebinom{n}{3} + \ebinom{n}{4} + \mathrm{etc} = 2^n
\]
haben; diese Reihe selbst ist natürlich die entwickelte Formel $(1+1)^n$. Wir wollen nun $\varphi = 180^{\circ}$ nehmen und es wird diese Reihe entstehen
\[
	1 - \ebinom{n}{1} + \ebinom{n}{2} - \ebinom{n}{3} + \ebinom{n}{4} - \mathrm{etc} = 0,
\]
diese Reihe ist natürlich $(1-1)^n$; Es sei auch $\varphi = 90^{\circ}$ und die daher entstehende Reihe wird
\[
	1 - \ebinom{n}{2} + \ebinom{n}{4} - \ebinom{n}{6} + \ebinom{n}{8} - \ebinom{n}{10} + \mathrm{etc}
\]
sein, deren Summe also
\[
	2^n(\cos{45^{\circ}})^n\cos{n45^{\circ}} = 2^{\frac{1}{2}n}\cos{n45^{\circ}}
\]
sein wird.
\paragraph{§15}
Diese letzte Reihe scheint umso größerer Aufmerksamkeit würdig, weil deren Gültigkeit nicht gerade wenig mysteriös ist; daher wird es nicht unpassend sein, dass einige Spezialfälle betrachtet werden und zwar für ganze positive Zahlen:
\begin{enumerate}
	\item Wenn $n=0$ ist, wird $1=1$ sein.
	\item Wenn $n=1$ ist, wird $1=\cos{45^{\circ}\sqrt{2}}$ sein.
	\item Wenn $n=2$ ist, wird $1-1=2\cos{90^{\circ}} = 0$ sein.
	\item Wenn $n=3$ ist, wird $1-3=2^{\frac{3}{2}}\cos{3\cdot 45^{\circ}} = -2$ sein.
	\item Wenn $n=4$ ist, wird $1-6+1=2^2\cos{4\cdot 45^{\circ}} = -4$ sein.
	\item Wenn $n=5$ ist, wird $1-10+5=2^{\frac{5}{2}}\cos{5\cdot 45^{\circ}} = -4$ sein.
	\item Wenn $n=6$ ist, wird $1-15+15-1=2^3\cos{6\cdot 45^{\circ}} = 0$ sein.
	\item Wenn $n=7$ ist, wird $1-21+35-7 = 2^{\frac{7}{2}}\cos{7\cdot 45^{\circ}} = 2^3$ sein.
	\item Wenn $n=8$ ist, wird $1-28+70+28+1 = 2^4\cos{8\cdot 45^{\circ}} = 2^4$ sein.
\end{enumerate}
etc.
\paragraph{§16}
Größere Aufmerksamkeit verdienen die Fälle, in denen für $n$ eine negative Zahl angenommen wird, in denen ja unendliche Reihen hervorgehen:
\begin{enumerate}
	\item Wenn $n=-1$ ist, wird
	\[
		1-1+1-1+1-1+1-1+1-1+\mathrm{etc} = \frac{\cos{45^{\circ}}}{\sqrt{2}} = \frac{1}{2}
	\]
	sein.
	\item Wenn $n=-2$ ist, wird
	\[
		1-3+5-7+9-11+14-15+17-\mathrm{etc} = \frac{\cos{2\cdot 45^{\circ}}}{2} = 0
	\]
	sein.
	\item Wenn $n=-3$ ist, wird
	\[
		1-6+15-28+45-66+91-\mathrm{etc} = \frac{\cos{3\cdot 45^{\circ}}}{\sqrt{8}} = -\frac{1}{4}
	\]
	sein.
	\item Wenn $n=-4$ ist, wird
	\[
		1-10+35-84+165-286+455-\mathrm{etc} = \frac{\cos{4\cdot 45^{\circ}}}{4} = -\frac{1}{4}
	\]
	sein.
	\item Wenn $n=-5$ ist, wird 
	\[
		1-15+70-210+495-1001+\mathrm{etc} = \frac{\cos{5\cdot 45^{\circ}}}{2^{\frac{5}{2}}} = -\frac{1}{8}
	\]
	sein.
	\item Wenn $n=-6$ ist, wird
	\[
		1-21+126-462+1287-3003+\mathrm{etc} = \frac{\cos{6\cdot 45^{\circ}}}{8} = 0
	\]
	sein, etc.
\end{enumerate}
Im Allgemeinen wird aber für diese Fälle
\begin{align*}
	&1 - \frac{\lambda (\lambda +1)}{1\cdot 2} + \frac{\lambda (\lambda +1)(\lambda +2)(\lambda +3)}{1\cdot 2\cdot 3\cdot 4} - \frac{\lambda (\lambda +1)(\lambda +2)(\lambda +3)(\lambda +4)(\lambda +5)}{1\cdot 2\cdot 3\cdot 4\cdot 5\cdot 6} \\
	&+\frac{\lambda (\lambda +1)\cdots (\lambda +7)}{1\cdot 2\cdots 8} - \frac{\lambda (\lambda +1)\cdots (\lambda +9)}{1\cdot 2\cdots 10} + \mathrm{etc}
\end{align*}
sein, die Summe welcher Reihe also gleich
\[
	\frac{\cos{\lambda \cdot 45^{\circ}}}{2^{\frac{1}{2}\lambda}}
\]
sein wird.
\section*{Problem 2}
Nachdem diese Reihe der Sinus vorgelegt wurde:
\[
	\frac{n}{1}\sin{\varphi} + \frac{n}{1}\cdot\frac{n-1}{2}\sin{2\varphi} + \frac{n}{1}\cdot\frac{n-1}{2}\cdot\frac{n-2}{3}\sin{3\varphi} + \mathrm{etc} = s,
\]
sodass durch die oben verwendeten Charaktere
\[
	s = \ebinom{n}{0}\sin{0\varphi} + \ebinom{n}{1}\sin{1\varphi} + \ebinom{n}{2}\sin{2\varphi} + \ebinom{n}{3}\sin{3\varphi} + \mathrm{etc}
\]
ist, ist der Wert dieser Summe $s$ reell auszudrücken.
\section*{Lösung}
\paragraph{§17}
Wenn wir also hier wiederum die Buchstaben
\[
	p = \cos{\varphi} + \sqrt{-1}\sin{\varphi}
\]
und
\[
	q = \cos{\varphi}-\sqrt{-1}\sin{\varphi}
\]
einführen, teile man, weil ja
\[
	p^n-q^n = 2\sqrt{-1}\sin{n\varphi}
\]
ist, die vorgelegte Reihe in die $2$ folgenden auf
\[
	2s\sqrt{-1} = \left. \begin{cases}
	\quad +\ebinom{n}{1}p + \ebinom{n}{2}pp + \ebinom{n}{3}p^3 + \ebinom{n}{4}p^4  \\
	\quad -\ebinom{n}{1}q - \ebinom{n}{2}qq - \ebinom{n}{3}q^3 - \ebinom{n}{4}q^4
	\end{cases}\right\}
	\mathrm{etc},
\]
woher natürlich
\[
	2s\sqrt{-1} = (1+p)^n - (1+q)^n
\]
sein wird.
\paragraph{§18}
Hier wird es wiederum gleich förderlich sein bemerkt zu haben, dass
\[
	1+p = (\sqrt{p} + \sqrt{q})\sqrt{p}
\]
und
\[
	1+q = (\sqrt{p} + \sqrt{q})\sqrt{q}
\]
ist, nach Verwendung welcher Werte
\[
	2s\sqrt{-1} = (\sqrt{p} + \sqrt{q})^n(p^{\frac{n}{2}} - q^{\frac{n}{2}})
\]
sein wird. Weil daher
\[
	p^{\frac{n}{2}} - q^{\frac{n}{2}} = 2\sqrt{-1}\sin{\tfrac{1}{2}n\varphi}
\]
ist und
\[
	\sqrt{p} + \sqrt{q} = 2\cos{\tfrac{1}{2}\varphi},
\]
wird daher, indem man durch $2\sqrt{-1}$ teilt, die gesuchte Summe reell ausgedrückt hervorgehen:
\[
	s = 2^n\cos^n{\tfrac{1}{2}\varphi}\sin{\tfrac{1}{2}n\varphi}.
\]
\end{document}